\documentclass{article}
\usepackage{latexsym}
\usepackage[dvips]{graphicx}

\begin{document}
\renewcommand{\thefootnote}{}
\footnotetext{
msc: transformation groups, manifolds\\
\hspace*{4.5mm}
keywords: orbits, Lie group actions, weakly embedded manifolds}
\newtheorem{Cor}{Corollary}

\newenvironment{NoNo} [1]
         {\medskip\par{\bf #1 }\it}{\rm \medskip \newline}

\newenvironment {prqed}        
         {\begin{description}\item[Proof:]}
         {\newline\nopagebreak\hfill $\Box$ \end{description}\medskip}

\renewcommand   {\phi}  {\varphi}
\newcommand     {\R}    {{\ifmmode{\rm I}\mkern-4mu{\rm R}
    \else\leavevmode\hbox{I}\kern-.17em\hbox{R}\fi}}
\newcommand     {\C}    {\ifmmode{{\rm C}\mkern-15mu{\phantom{\rm t}\vrule}}
    \mkern10mu
    \else\leavevmode\hbox{C}\kern-.5em\hbox{I}\kern.3em\fi}
\newcommand     {\N}    {\ifmmode{\rm I}\mkern-3.5mu{\rm N}
    \else\leavevmode\hbox{I}\kern-.16em \hbox{N}\fi}
\newcommand     {\Z}    {\Bbb{Z}}
\newcommand     {\Q} {{\cal Q}}

\hyphenation{}

\begin{center}

{\Large
Orbits of Lie Group Actions are Weakly Embedded
}\bigskip\\
Domenico P.L. Castrigiano\\

Sandra A. Hayes\\
Zentrum Mathematik der Technischen Universit\"at M\"unchen\\
Arcisstr. 21, 80333 Munich, Germany\smallskip\\
{\tt
castrig@mathematik.tu-muenchen.de\\
hayes@mathematik.tu-muenchen.de}\\

\end{center}
\bigskip
\begin{quote}
{\bf Abstract.}
In this note we prove that whenever a Lie group $G$ acts on a manifold $X$, then the orbit
$Gx$ through any point $x$ of $X$ is a weakly embedded submanifold
of $X$. The investigation of this  problem was inspired by an application to Catastrophe Theory.
\end{quote}

%
\section{Introduction}
When a Lie group $G$ acts on a manifold $X$, the orbit $Gx$ through a point $x$ of $X$ is a
special immersed submanifold of $X$, namely it is weakly embedded. According 
to \cite{P}, an immersed
submanifold $Y$ of $X$ is said to be 
{\it weakly embedded, if every smooth map $g:S\to X$ of a manifold $S$ into $X$ with image in $Y$ induces
  a smooth map $g:S\to Y$.}
In particular, a  smooth curve $c$ in $X$, whose image lies in $Y$, is actually
 a smooth curve in $Y$. Immersed submanifolds, which are not reglar submanifolds, i.e. whose topology is strictly finer than the induced one, need not possess this property. A simple illustration is the following Figure Eight.
\bigskip\bigskip

$$
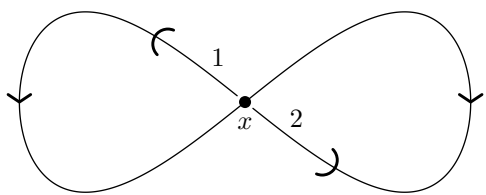
\includegraphics{eight.1}$$
\bigskip

The curve, whose image is composed of parts 1 and 2 as well as the point $x$,
lies in the Figure Eight but is not continuous as a map into the Figure Eight \cite[1.31]{W}.
Consequently, this Figure Eight is not weakly embedded.

The main step in proving that orbits of Lie group actions are weakly embedded 
consists in showing 
that any  $C^1$ curve whose image lies in an orbit locally coincides with the image of a curve
in the group under the orbital map (Lifting Theorem). The  idea for studying this problem
arose from analyzing its relevance to one of the major theorems in Catastrophe Theory,
namely John Mather's \cite{M} necessary condition for determining when a smooth function
of several variables coincides with its own Taylor polynomial up to a local diffeomorphism
(see \cite[Chap. 4]{Ca/Ha}).
%
\section{The Linearization Lemma}

In this note a {\it manifold} $X$ means a smooth manifold, i.e. a $C^\infty$ manifold.
$X$ need not be Hausdorff nor have a countable base. An {\it immersed submanifold}
$Y$ of $X$ is a manifold which is a subset of $X$ such that the inclusion $i:Y\to X$
is an immersion. Clearly, the topology of $Y$ is finer than the relative topology of $X$.
When these topologies coincide, $Y$ is called a
{\it regular submanifold}.
Weakly embedded submanifolds are intermediary between regular and immersed submanifolds
(see \cite{Mo}, \cite{S}). Note that in the definition of weakly embedded, the map
$g:S\to X$ is smooth if it is continuous \cite[1.13]{W}.
More generally, for every $k$ it follows as in \cite{W} that every $C^k$ map $g:S\to X$
with image in $Y$ induces a $C^k$ map $g:S\to Y$ if $g:S\to Y$ is only continuous.
Furthermore, to notify that an immersed submanifold $Y$ of $X$ is weakly embedded, i.e.
possesses the lifting property of smooth maps given above, it suffices to consider just the
lifting of $C^1$ curves:

\begin{NoNo}{Proposition}
Let $Y$ be an immersed submanifold of a manifold $X$ such that every $C^1$ curve in $X$ with
image in $Y$ induces a continuous curve in $Y$. Then $Y$ is weakly embedded in $X$.
\end{NoNo}

\begin{prqed}
Otherwise, for a manifold $S$ and a smooth map $g:S\to X$ with $g(S)\subset Y$, there would
be a sequence $(s_n)$ in $S$ converging to a point $s_0$ in $S$ whereas its image
$\bigl(g(s_n)\bigr)$ would not converge to $g(s_0)$ in $Y$. 
In charts, let $s_0$ be the origin. Without restriction,
$t_n:=|s_n|<1$ is strictly monotonically decreasing.
By taking chords between successive points $s_n$ and $s_{n+1}$ and smoothing
the corners, one obtains a smooth curve $c:\;]0,1[\to S$ with $c(t_n)=s_n$ which is
continuously extendable to the origin by setting $c(0):=0$. A further refinement involving a
reparametrisation using the third power of the arc length results in a smooth curve
$c:\;]0,1[\to S$ with $c(t_n)=s_n$ which is $C^1$ extendable to $]-1,1[$ by defining $c(t)=0$,
$t\leq0$. Thus, $g\circ c:\;]-1,1[\to X$ is $C^1$ with image in $Y$ but is not continuous as a
map into $Y$, contradicting the hypothesis of the proposition.
\end{prqed}

Leaves of a foliation are classical examples of weakly embedded submanifolds
(see \cite{Mo}, \cite{S}). The objective here is to prove that if a Lie group $G$ acts on a
manifold $X$, the orbits $Gx$ through points $x$ in $X$ are weakly embedded. This will be
proved by showing that the property of $C^1$ curves given in the above proposition is
satisfied. The proof of this requires a Linearization Lemma and a Lifting Theorem.
The latter will be proved by solving a differential equation. Tangent vector will always refer to  vectors tangent to  $C^1$ curves.

A Lie group $G$ is assumed to be separable and, consequently, $\sigma$--compact.
Let $(\xi,x) \mapsto \xi x$ be a smooth action of $G$ on $X$. Recall that for $x \in X$
the canonical map $G / G_x \to X, \; \xi G_x \mapsto \xi x$, is an injective immersion and
its image $Gx$ is an immersed submanifold of $X$ \cite[III,\thinspace 1.7]{B}.
Denote by $\rho : G \to X$ the orbital map $\xi \mapsto \xi x$. The rank $r$ of $\rho$ is
constant \cite[III,\thinspace 1.5]{B}. Therefore $T_{\xi x}(Gx)=d_{\xi}\rho(T_{\xi}G)$ implies
\begin{equation}\label{eq1}
r = \mbox{dim}\: T_{\xi x}(Gx)
\end{equation}
independent of $\xi \in G$.
%
%
\begin{NoNo}{Linearization Lemma}
Any $C^k$ curve $c$ in $X$ based at a point $x=c(0)$ in $X$, whose image lies in the orbit
$Gx$ through $x$, is tangent at $x$ to a $C^k$ curve in that orbit, i.e. $d_0c$ is an element
in the tangent space $T_x(Gx)$ of the orbit at $x$, for $k\geq1$.
\end{NoNo}
\begin{prqed}
Let $d$ be the dimension of $G$ and $N$ that of $X$ around $x \in X$. By (\ref{eq1}),
the Rank Theorem can be applied to the orbital map $\rho$ so that in charts $u$ for $G$ at
$e$ with $u(e)=0$ and $v$ for $X$ at $x$ with $v(x)=0$, $\rho$ becomes $p \mapsto (pr_r(p),0)$
for $0 \in \R^{N-r}$, where $pr_r$ is the projection onto the first $r$ coordinates.
Let $c$ be defined on $I=]-1,1[$.

For any $\epsilon \in ]0,1[$ there exist an $s \in I$ with $|s| \leq \epsilon$ and convergent
sequences $(t_n)$ in $I \backslash \{ s \}$ and $(\xi_n)$ in $G$ such that $t_n \to s$
and $c(t_n) = \rho(\xi_n)$ hold. To see this, for every $t$ in
$I_{\epsilon} := [-\epsilon,\epsilon]$ choose any $\xi_t$ in $G$ with $c(t)=\rho(\xi_t)$.
The $\sigma$--compactness of $G$ insures the existence of infinitely many $t \in I_{\epsilon}$
such that all $\xi_t$ lie in some compact subset of $G$. Choose a sequence $(t_l)$ of such
$t$'s for which $(\xi_{t_l})$ converges. Then any convergent subsequence of $(t_l)$ proves
the claim.

Assume first that $s=0$. Let $\xi := \lim_n \xi_n$ and put $\xi'_n := \xi_n \xi^{-1}$.
Then $\xi'_n \to e$ and $c(t_n) = \rho(\xi'_n)$ holds. Therefore
$v \circ c \: (t_n) \in \R^r \times \{ 0 \}$ and $(v \circ c)'(0) = (k,0)$
for some $k \in \R^r$ and $0 \in \R^{N-r}$ due to $v \circ c (0) = 0$.
Now $C(t) := u^{-1}(tk,0),\: 0 \in \R^{d-r}$, defines a curve $\rho \circ C$ in the orbit
$Gx$ satisfying $v \circ \rho \circ C(t) = (tk,0)$. Since it also has $(k,0)$ as its tangent 
at $x$, we are finished.

Next consider $s\not=0$. Let $\rho_t : G \to X$ denote the orbital map $\xi \mapsto \xi c(t)$ of
$G$ into $X$ defined by $c(t)$ for $t$ in $I$. Obviously, $\rho_0 = \rho$, and by hypothesis
$Gc(t)=Gx$. Since $c(t_n) = \rho_s(\xi'_n)$, similar to the argument for $s=0$ we obtain
\begin{equation}\label{eq2}
d_s c \in T_{c(s)} (Gx)\:.
\end{equation}

An appropriate continuity argument will now show that (\ref{eq2}) holds for $s=0$.
For every $t$ in $I$ the dimension of the tangent space $T_t := T_{c(t)}(Gx) = d_e\rho_t(T_eG)$
is $r$ by (\ref{eq1}). Choose elements $\Gamma_1,\dots,\Gamma_r$ in $T_eG$ so that
$d_e\rho(\Gamma_j), \: 1 \leq j \leq r$, form a basis of $T_0$. Since the maps
\[
I \to TX, \: t \mapsto \Gamma_j(t) := d_e\rho_t(\Gamma_j)
\]
are continuous, the $\Gamma_j(t), \: 1 \leq j \leq r$,
constitute a basis of $T_t$ for small $t$. The following calculations are best done in charts.
Using an orthogonal basis obtained from $\Gamma_j(t), \: 1 \leq j \leq r$, an expression for
the length $l(t)$ of the component of $d_tc$ in $\R^N$ orthogonal to $T_t$ is deduced which
depends continuously on $t$. For the special parameters $s$, (\ref{eq2}) is equivalent to
$l(s)=0$, and we obtain $l(0)=0$, implying $d_0c \in T_0$.
\end{prqed}

A special case of this Linearization Lemma is treated in [Ca/Ha, 4.37].

%
\section{The Lifting Theorem}
Even more than is claimed in the Linearization Lemma holds. Around the origin,
$c$ actually coincides with  $\rho \circ C$ where $C$ is a curve in $G$.

\begin{NoNo}{Lifting Theorem}
Let $G$ be a Lie group acting on a manifold $X$. Then any $C^k$ curve $c$ in $X$ based at a
point $x$ in $X$, whose image lies in the orbit $Gx$ through $x$, can be locally lifted to a
$C^k$ curve $C$ in $G$ based at the identity, i.e., $C(t)x=c(t)$ for small $t$, for $k\geq1$.
\end{NoNo}
\begin{prqed}
The curve $C$ will be the solution of a differential equation. To set up this equation choose
a right inverse $r(t) : T_t \to T_eG$ of $d_e\rho_t$, i.e. $(d_e\rho_t)r(t)\gamma = \gamma$
for all $\gamma \in T_t$, which depends smoothly on $t$ for small $t$. Let
$R_{\xi} : G \to G, \: \eta \mapsto \xi\eta$, denote right multiplication on $G$.
Because the Linearization Lemma implies $c'(t) := d_tc \in T_t$, the following differential
equation for a curve $C$ in $G$ based at $e$ is well-defined:
\begin{equation}\label{eq3}
C'(t) = d_e R_{C(t)} r(t) c'(t) \: .
\end{equation}
It remains to check that $C(t)c(0) = c(t)$ holds or equivalently $\bigl(C(t)^{-1}c(t)\bigr)'=0$
for small $t$. For $\xi \in G$ let $\alpha_{\xi}$ denote the diffeomorphism
$X \to X, \: y \mapsto \xi y$. Then
\begin{equation}\label{eq4}
\bigl( C(t)^{-1} c(t)\bigr)' 
  = d_{C(t)^{-1}}\:\rho_t\bigl(C(t)^{-1}\bigr)'+d_{c(t)}\:\alpha_{C(t)^{-1}}\:c'(t)\:.
\end{equation}
The second summand is the negative of the first. To see this, observe that (\ref{eq3}) implies
\[
c'(t) = d_e \: \rho_t (d_e R_{C(t)})^{-1} \: C'(t) \: .
\]
Differentiating $e = C(t) C(t)^{-1}$ yields
$0 = ( d_e R_{C(t)} )^{-1} C'(t) + d_{C(t)^{-1}} L_{C(t)} (C(t)^{-1})' \: ,$
where $L_{\xi}$ denotes left multiplication on $G$ and where
$d_{C(t)} R_{C(t)^{-1}} = (d_e R_{C(t)})^{-1}$ has been used.
The second summand in (\ref{eq4}) now becomes
\[
-d_{c(t)} \: \alpha_{C(t)^{-1}} \: d_e \: \rho_t \: d_{C(t)^{-1}} \: L_{C(t)} (C(t)^{-1})' \: .
\]
However, since $\rho_t = \alpha_{C(t)^{-1}} \circ \rho_t \circ L_{C(t)}$ is true,
this term is the negative of the first summand in (\ref{eq4}).
\end{prqed}
A special case of this theorem is in [Ca/Ha, Exercise 4.7].
Because the orbital map $G \to Gx, \: \xi \mapsto \xi x$, is smooth, we can conclude with
\begin{Cor}
Every $C^k$ curve in $X$ based at a point $x$ in $X$, whose image lies in the orbit $Gx$
through $x$, is a $C^k$ curve in that orbit, for $k\geq1$.
\end{Cor}
\begin{Cor}
Every $C^k$ curve in an immersed submanifold $Z$ of $X$, which is contained in the orbit
$Gx$, defines a $C^k$ curve in that orbit, for $k\geq1$.
\end{Cor}
\begin{Cor}
The orbits of $G$ in $X$ are weakly embedded.
\end{Cor}

\end{document}